\documentclass{article}

\usepackage[a4paper, left=3.3cm,top=3.3cm,right=3.3cm,bottom=3.7cm]{geometry}
\usepackage{concmath}
\usepackage[T1]{fontenc}

\usepackage{amsmath,amsfonts,amssymb}
\usepackage{graphicx}
\usepackage[colorlinks=true, allcolors=blue]{hyperref}
\usepackage{mathtools}
\usepackage{framed}
\usepackage{enumitem}
\setenumerate{topsep=0em}
\setitemize{topsep=0em,label=\scriptsize{$\blacksquare$}}

\usepackage{siunitx}
\usepackage{authblk}
\usepackage{subcaption}
\usepackage{soul}

\usepackage{tikz}
\usepackage{tikz-3dplot}
\usetikzlibrary{external}
\usepackage{booktabs}

\newcommand*{\supp}{\mathrm{supp}}

\newcommand*{\R}{\mathbb{R}}

\newcommand{\ro}{\mathbf{r}}
\newcommand{\soo}{\mathbf{s}}

\newcommand{\A}{\mathcal{A}}

\newcommand{\z}{\mathbf{z}}
\newcommand{\y}{\mathbf{y}}
\newcommand{\x}{\mathbf{x}}

\newcommand{\f}{\mathbf{f}}

\newcommand*{\Int}[4]{\int_{#1}^{#2}\!{#3}\,\mathrm{d}{#4}}

\renewcommand{\phi}{\varphi}
\renewcommand{\epsilon}{\varepsilon}

\DeclarePairedDelimiter{\abs}{\lvert}{\rvert}
\DeclarePairedDelimiter{\norm}{\lVert}{\rVert}

\DeclareMathOperator*{\argmin}{argmin}

\DeclareMathOperator{\prox}{prox}

\DeclareMathOperator{\sign}{sign}

\newcommand{\edot}{\,\cdot\,}

\definecolor{drot}{rgb}{0,0,0}
\definecolor{blor}{rgb}{1,0,1}
\definecolor{blgr}{rgb}{0,1,1}
\definecolor{bblau}{rgb}{0,0,1}
\definecolor{goyel}{rgb}{0.3,0,1}
\definecolor{orred}{rgb}{0,0.39,0}
\newcommand{\hlc}[2][gray]{ {\sethlcolor{#1} \hl{#2}} }

\def\Put(#1,#2)#3{\leavemode\makebox(0,0){\put(#1,#2){#3}}}

\def\mirror at (#1,#2){\filldraw[color=orange](#1,#2-0.05)--(#1+0.15,#2-0.05)--(#1+0.15,#2-0.1)--(#1+0.25,#2)--(#1+0.15,#2+0.1)--(#1+0.15,#2+0.05)--(#1,#2+0.05)--cycle;}
\def\mirrordown at (#1,#2){\filldraw[color=lime](#1-0.02,#2-0.02)--(#1+0.12,#2-0.16)--(#1+0.085,#2-0.195)--(#1+0.21,#2-0.21)--(#1+0.195,#2-0.085)--(#1+0.16,#2-0.12)--(#1+0.02,#2+0.02)--cycle;}
\def\mirrord at (#1,#2){\filldraw[color=lime](#1-0.05,#2)--(#1-0.05,#2-0.15)--(#1-0.1,#2-0.15)--(#1,#2-0.25)--(#1+0.1,#2-0.15)--(#1+0.05,#2-0.15)--(#1+0.05,#2)--cycle;}
\def\mirroru at (#1,#2){\filldraw[color=blgr](#1-0.05,#2)--(#1-0.05,#2+0.15)--(#1-0.1,#2+0.15)--(#1,#2+0.25)--(#1+0.1,#2+0.15)--(#1+0.05,#2+0.15)--(#1+0.05,#2)--cycle;}
\def\mirrortwo at (#1,#2){\filldraw[color=orred](#1,#2-0.05)--(#1+0.15,#2-0.05)--(#1+0.15,#2-0.1)--(#1+0.25,#2)--(#1+0.15,#2+0.1)--(#1+0.15,#2+0.05)--(#1,#2+0.05)--cycle;}
\def\mirrorthree at (#1,#2){\filldraw[color=blgr](#1-0.1,#2-0.05)--(#1+0.15,#2-0.05)--(#1+0.15,#2-0.1)--(#1+0.25,#2)--(#1+0.15,#2+0.1)--(#1+0.15,#2+0.05)--(#1-0.1,#2+0.05)--cycle;}
\def\mirrorfour at (#1,#2){\filldraw[color=gray](#1-0.1,#2-0.05)--(#1+0.15,#2-0.05)--(#1+0.15,#2-0.1)--(#1+0.25,#2)--(#1+0.15,#2+0.1)--(#1+0.15,#2+0.05)--(#1-0.1,#2+0.05)--cycle;}
\def\mirrorfive at (#1,#2){\filldraw[color=lime](#1-0.1,#2-0.05)--(#1+0.15,#2-0.05)--(#1+0.15,#2-0.1)--(#1+0.25,#2)--(#1+0.15,#2+0.1)--(#1+0.15,#2+0.05)--(#1-0.1,#2+0.05)--cycle;}

\allowdisplaybreaks

\title{NETT Regularization for Compressed Sensing Photoacoustic Tomography}

\author{Stephan Antholzer}
\author{Johannens Schwab}

\affil{Department of Mathematics, University of Innsbruck\authorcr
Technikerstrasse 13, 6020 Innsbruck, Austria\authorcr
E-mail: {\tt \{stephan.antholzer,johannes.schwab@uibk.ac.at\}@uibk.ac.at}}

\author{Johannes Bauer-Marschallinger}
\author{Peter Burgholzer}

\affil{Department of Mathematics, University of Innsbruck\authorcr
Technikerstrasse 13, 6020 Innsbruck, Austria\authorcr
E-mail: {\tt \{j.bauer-marschallinger,peter.burgholzer\}@recendt.at}}

\author{Markus Haltmeier}

\affil{Department of Mathematics, University of Innsbruck\authorcr
Technikerstrasse 13, 6020 Innsbruck, Austria\authorcr
E-mail: {\tt markus.haltmeier@uibk.ac.at}}

\begin{document}
\maketitle

\begin{abstract}
We discuss several methods for image reconstruction in  compressed sensing  photoacoustic tomography (CS-PAT). In particular, we apply the  deep learning method  of  [H. Li, J. Schwab, S. Antholzer, and M. Haltmeier. \emph{NETT: Solving Inverse Problems with Deep Neural Networks (2018)}, arXiv:1803.00092], which is based on a learned regularizer,   for the first time  to the  CS-PAT problem.  We propose a network architecture and training strategy
for the NETT that we expect to be useful for other inverse problems as well.
All algorithms are compared and evaluated on simulated data, and  validated using experimental data for two
different types of phantoms. The results on the one the hand indicate great potential of deep learning methods,
and on the other hand show that significant future work is required to improve their performance
on real-word data.

\bigskip\noindent
\textbf{Keywords:}
Compressed sensing, photoacoustic tomography, deep learning,  NETT, learned regularizer, Tikhonov regularization, $\ell^1$-regularization, neural networks, inverse problems

\end{abstract}

\section{Introduction}

Compressed Sensing (CS) is a promising tool to reduce the number spatial measurements in photoacoustic tomography (PAT), while still keeping good image quality. Reducing the number of measurements can be used to lower system costs, to speed up data acquisition, and to reduce motion artifacts \cite{sandbichler2015novel,arridge2016accelerated,haltmeier2016compressed}.
In this work, we concentrate on 2D PAT, which arises in PAT with integration line detectors \cite{burgholzer2007temporal,paltauf2007photacoustic}. If one can use enough sensor locations  such that Shannon's sampling theory is applicable, it is possible to get high resolution  images using analytic inversion methods \cite{haltmeier2016sampling}. In practice, however,  the number of measurements is much lower than required for high resolution images according to
Shannon's sampling theory. In such a situation, standard reconstruction methods result in
undersampling artifacts and low resolution images.

To accelerate the data acquisition process, systems with 64 line detectors have been built~\cite{gratt201564line,beuermarschallinger2015photacoustic}.
Such systems offer the possibility to collect 2D photoacoustic (PA) projection images  of the 3D source images at a frame-rate  of  $\SI{20}{\hertz}$ or higher \cite{paltauf2017piezoelectric}. 
Using 64 spatial sampling positions still results in
highly under-sampled data. In order to get high resolution reconstructions
from such data, one has to exploit additional information available on the PA
source images. In~\cite{sandbichler2018sparsification}, the sparsity of the Laplacian is used in 
combination with $\ell_1$-minimization.
Recently, machine learning methods have  been applied to CS-PAT \cite{japan2018proc}.
In this work, we develop several machine learning methods for CS-PAT
and apply them  to experimental CS-PAT data. Comparison with 
joint $\ell_1$-minimization is also given.

In contrast to \cite{japan2018proc}, we also implement the NETT
\cite{li2018nett}, which uses a neural network as trained  regularizer for CS-PAT image reconstruction.
In particular, we propose a simpler network architecture and training strategy
than the one used  in \cite{li2018nett}. The proposed strategy for the
regularizer in NETT may be useful for other inverse problems as well.

\section{Background}
\label{sec:background}

In this section  we  describe the CS-PAT problem, and present the 
joint $\ell^1$-algorithm as well as 
the standard deep learning approach for CS-PAT image reconstruction.
The NETT will be introduced  in Section~\ref{sec:nett}.

\subsection{Compressed sensing PAT}

PAT relies on the following principle. When an object of interest is illuminated with a short  laser pulse,  an acoustic pressure wave is
induced inside the object.
This pressure wave is recorded on the outside of the object and used for reconstructing an image of the interior. In the following we will restrict ourselves to the 2D case in  a circular measurement geometry, which arises when using integrating line detectors~\cite{burgholzer2007temporal,paltauf2007photacoustic}.

Let $p_0\colon\R^2\to\R$ denote  the PA source (initial pressure distribution). The induced  pressure wave 
satisfies the following equation
\begin{equation}\label{eq:wave}
    \partial^2p(\ro, t) - c^2 \Delta_\ro p(\ro, t)
    = \delta^\prime(t) p_0(\ro) \quad \text{for } (\ro, t)\in \R^2\times \R_+\,,
\end{equation}
where $\ro\in\R^2$ is the spatial location, $t\in\R_+$ the time variable, $\Delta_\ro$ the spatial Laplacian, and $c$ is the  speed of sound. We further assume that the PA source $p_0(\ro)$ vanishes outside the disc $B_R = \{ \x \in \R^2 \mid \norm{\x}<R\}$ and we set $p(\ro,t)=0$ for $t<0$.
Then $p(\ro,t)$ is uniquely defined and we refer to it as the causal solution of~\eqref{eq:wave}.

The PAT problem consists in recovering $p_0(\ro)$ from measurements of $p(\soo,t)$ on $\partial B_R \times (0, \infty)$, where $\soo$ stands for the  detector location. 
In the full data case, as shown in~\cite{FinHalRak07}, the following 
filtered backprojection formula (FBP) formula yields an exact reconstruction of the PA source,
\begin{equation}\label{eq:fbp}
    p_0(\ro) = -\frac{1}{\pi R} \Int{\partial B_R}{}{\Int{\abs{\ro-z}}{\infty}{\frac{(\partial_ttp)(\soo,t)}{\sqrt{t^2-\abs{\ro-\soo}^2}}}{t}}{S(\soo)} \,.
\end{equation}
Additionally, in \cite{FinHalRak07}  it   was shown that the  operator defined by the right hand side of  \eqref{eq:fbp} is the adjoint  of the forward operator of the wave equation.

In practical application, the pressure can only be measured with a finite number of samples. This means we measure data
\begin{equation}\label{eq:mess1}
    p(\soo_k,t_l) \text{ for } (k,l)\in\{1,\ldots,M\}\times \{1,\ldots,Q\} \,,
\end{equation}
where the sampling  points are uniformly sampled, i.e.
\begin{equation} \label{eq:mess2}
    \soo_k = \left[\begin{array}{c} R\cos(2\pi(k-1)/M)\\ R\sin(2\pi(k-1)/M) \end{array}\right]
	\quad t_l = 2R(l-1)/(Q-1)\,.
\end{equation}
As shown in \cite{haltmeier2016sampling},
using classical Shannon  sampling theory,
the number $M$ of spatial measurements   in  \eqref{eq:mess1}, \eqref{eq:mess2}  determines the resolution of the reconstructed PA source.
To reduce the number of measurements while preserving high resolution  we apply CS. Instead of collecting $M$ samples, we
measure generalized samples
\begin{equation}\label{eq:CS}
    \mathbf{g}(j,l) = \sum_{k=1}^M \mathbf{S}[j,k]p(\ro_k, t_l) \quad \text{for } j\in\{1,\ldots,m\}\,,
\end{equation}
with $m\ll M$. Several choices of the sampling matrix $\mathbf{S}$ exist~\cite{sandbichler2015novel,arridge2016accelerated,haltmeier2016compressed}. In this work we will focus on two cases, namely a
deterministic subsampling matrix and a random Bernoulli matrix.

Let us denote the discretized solution operator of the wave equation by $\mathcal{W}\in\R^{MQ\times n}$, were $n$ is the discretization size of 
the reconstruction,  and by $\mathcal{S}=\mathbf{S}\otimes \mathbf{I}\in\R^{mQ\times MQ}$ the Kronecker product between the CS measurement matrix $\mathbf{S}$ and the identity matrix.
If we denote the discrete initial pressure by $\x\in\R^n$, we can write the measurement process the following way
\begin{equation}\label{eq:cspat}
    \y = \A\x + \varepsilon \quad\text{with } \A=\mathcal{S}\circ\mathcal{W}
    \in\R^{mQ\times n} \,.
\end{equation}
Since $mQ\ll n$, equation \eqref{eq:cspat} is  under-determined  and its solution requires specialized reconstruction algorithms that
incorporate additional knowledge about the unknowns.
Such algorithms will be described in the following.

\subsection{Residual networks}

Deep learning has been recently applied to several image reconstruction problems~\cite{chen2017lowdose,jin2017deep,han2016deep} including  PAT~\cite{japan2018proc,antholzer2018deep,antholzer2018photoacoustic,hauptmann2017model,waibel2018reconstruction}.

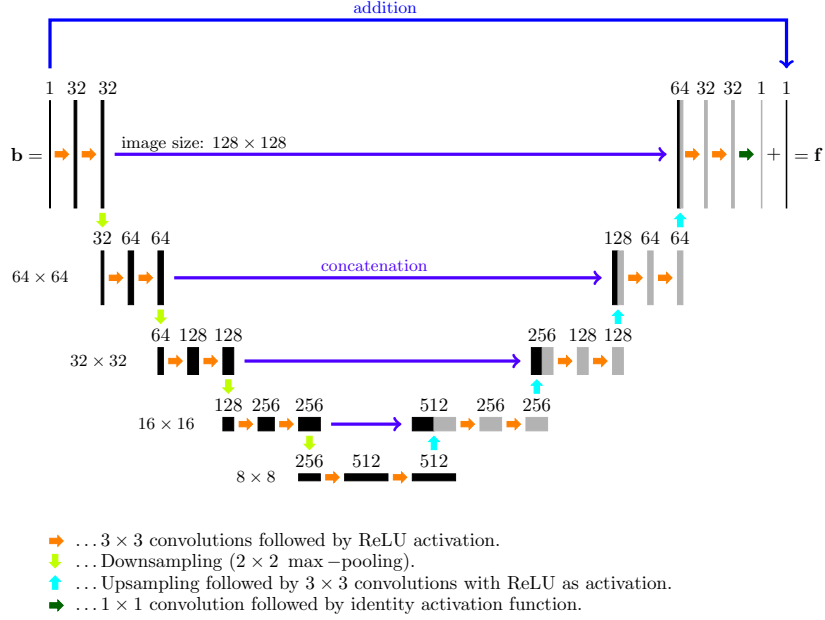
\begin{figure}[htp!]
    \centering
    \resizebox{!}{0.6\textwidth}{
	\begin{tikzpicture}
\draw(0,3)node[left]{};
\draw(0,-7)node[left]{};
\draw(1.2,0.2)node[right]{\small{\color{black}{image size: $128\times 128$}}};
\filldraw[color=black](0,1)--(0,-1)--(0.01,-1)--(0.01,1)--cycle;
\draw(0.005,1)node[above]{$1$};
\mirror at (0.1,0)
\filldraw[color=black](0.45,1)--(0.45,-1)--(0.5,-1)--(0.5,1)--cycle;
\draw(0.5,1)node[above]{$32$};
\mirror at (0.6,0)
\filldraw[color=black](0.95,1)--(0.95,-1)--(1,-1)--(1,1)--cycle;
\draw(1.075,1)node[above]{$32$};
\mirrord at (0.975,-1.1)
\draw(0.475,-2.3)node[left]{\small{\color{drot}{$64\times 64$}}};
\filldraw[color=black](0.95,-1.8)--(0.95,-2.8)--(1,-2.8)--(1,-1.8)--cycle;
\draw(0.975,-1.8)node[above]{$32$};
\mirror at (1.1,-2.3)
\filldraw[color=black](1.45,-1.8)--(1.45,-2.8)--(1.55,-2.8)--(1.55,-1.8)--cycle;
\draw(1.5,-1.8)node[above]{$64$};
\mirror at (1.65,-2.3)
\filldraw[color=black](2,-1.8)--(2,-2.8)--(2.1,-2.8)--(2.1,-1.8)--cycle;
\draw(2.05,-1.8)node[above]{$64$};
\mirrord at (2.05,-2.9)
\draw(1.55,-3.85)node[left]{\small{\color{drot}{$32\times 32$}}};
\filldraw[color=black](2,-3.6)--(2,-4.1)--(2.1,-4.1)--(2.1,-3.6)--cycle;
\draw(2.05,-3.6)node[above]{$64$};
\mirror at (2.2,-3.85)
\filldraw[color=black](2.55,-3.6)--(2.55,-4.1)--(2.75,-4.1)--(2.75,-3.6)--cycle;
\draw(2.65,-3.6)node[above]{$128$};
\mirror at (2.85,-3.85)
\filldraw[color=black](3.2,-3.6)--(3.2,-4.1)--(3.4,-4.1)--(3.4,-3.6)--cycle;
\draw(3.3,-3.6)node[above]{$128$};
\mirrord at (3.3,-4.2)
\draw(2.8,-5.025)node[left]{\small{\color{drot}{$16\times 16$}}};
\filldraw[color=black](3.2,-4.9)--(3.2,-5.15)--(3.4,-5.15)--(3.4,-4.9)--cycle;
\draw(3.3,-4.9)node[above]{$128$};
\mirror at (3.5,-5.025)
\filldraw[color=black](3.85,-4.9)--(3.85,-5.15)--(4.15,-5.15)--(4.15,-4.9)--cycle;
\draw(4,-4.9)node[above]{$256$};
\mirror at (4.25,-5.025)
\filldraw[color=black](4.6,-4.9)--(4.6,-5.15)--(5.0,-5.15)--(5.0,-4.9)--cycle;
\draw(4.8,-4.9)node[above]{$256$};
\mirrord at (4.8,-5.25)
\draw(4.3,-6.015)node[left]{\small{\color{drot}{$8\times 8$}}};
\filldraw[color=black](4.6,-5.95)--(4.6,-6.08)--(5,-6.08)--(5,-5.95)--cycle;
\draw(4.8,-5.95)node[above]{$256$};
\mirror at (5.1,-6.015)
\filldraw[color=black](5.45,-5.95)--(5.45,-6.08)--(6.25,-6.08)--(6.25,-5.95)--cycle;
\draw(5.85,-5.95)node[above]{$512$};
\mirror at (6.35,-6.015)
\filldraw[color=black](6.7,-5.95)--(6.7,-6.08)--(7.5,-6.08)--(7.5,-5.95)--cycle;
\draw(7.1,-5.95)node[above]{$512$};
\mirroru at (7.1,-5.5)
\filldraw[color=black](6.7,-4.9)--(6.7,-5.15)--(7.1,-5.15)--(7.1,-4.9)--cycle;
\filldraw[color=black!30](7.1,-4.9)--(7.1,-5.15)--(7.5,-5.15)--(7.5,-4.9)--cycle;
\draw(7.1,-4.9)node[above]{$512$};
\mirror at (7.6,-5.025)
\filldraw[color=black!30](7.95,-4.9)--(7.95,-5.15)--(8.35,-5.15)--(8.35,-4.9)--cycle;
\draw(8.15,-4.9)node[above]{$256$};
\mirror at (8.45,-5.025)
\filldraw[color=black!30](8.8,-4.9)--(8.8,-5.15)--(9.2,-5.15)--(9.2,-4.9)--cycle;
\draw(9,-4.9)node[above]{$256$};
\mirroru at (9,-4.45)
\filldraw[color=black](8.9,-3.6)--(8.9,-4.1)--(9.1,-4.1)--(9.1,-3.6)--cycle;
\filldraw[color=black!30](9.1,-3.6)--(9.1,-4.1)--(9.3,-4.1)--(9.3,-3.6)--cycle;
\draw(9.1,-3.6)node[above]{$256$};
\mirror at (9.4,-3.85)
\filldraw[color=black!30](9.75,-3.6)--(9.75,-4.1)--(9.95,-4.1)--(9.95,-3.6)--cycle;
\draw(9.85,-3.6)node[above]{$128$};
\mirror at (10.05,-3.85)
\filldraw[color=black!30](10.4,-3.6)--(10.4,-4.1)--(10.6,-4.1)--(10.6,-3.6)--cycle;
\draw(10.5,-3.6)node[above]{$128$};
\mirroru at (10.5,-3.15)
\filldraw[color=black](10.4,-1.8)--(10.4,-2.8)--(10.5,-2.8)--(10.5,-1.8)--cycle;
\filldraw[color=black!30](10.5,-1.8)--(10.5,-2.8)--(10.6,-2.8)--(10.6,-1.8)--cycle;
\draw(10.5,-1.8)node[above]{$128$};
\mirror at (10.7,-2.3)
\filldraw[color=black!30](11.05,-1.8)--(11.05,-2.8)--(11.15,-2.8)--(11.15,-1.8)--cycle;
\draw(11.1,-1.8)node[above]{$64$};
\mirror at (11.25,-2.3)
\filldraw[color=black!30](11.6,-1.8)--(11.6,-2.8)--(11.7,-2.8)--(11.7,-1.8)--cycle;
\draw(11.65,-1.8)node[above]{$64$};
\mirroru at (11.65,-1.35)
\filldraw[color=black](11.6,1)--(11.6,-1)--(11.65,-1)--(11.65,1)--cycle;
\filldraw[color=black!30](11.65,1)--(11.65,-1)--(11.7,-1)--(11.7,1)--cycle;
\draw(11.65,1)node[above]{$64$};
\mirror at (11.75,0)
\filldraw[color=black!30](12.1,1)--(12.1,-1)--(12.15,-1)--(12.15,1)--cycle;
\draw(12.125,1)node[above]{$32$};
\mirror at (12.25,0)
\filldraw[color=black!30](12.6,1)--(12.6,-1)--(12.65,-1)--(12.65,1)--cycle;
\draw(12.625,1)node[above]{$32$};
\mirrortwo at (12.75,0)
\filldraw[color=black!30](13.15,1)--(13.15,-1)--(13.16,-1)--(13.16,1)--cycle;
\draw(13.155,1)node[above]{$1$};
\draw(13.12,0)node[right]{$+$};
\filldraw[color=black](13.61,1)--(13.61,-1)--(13.62,-1)--(13.62,1)--cycle;
\draw(13.615,1)node[above]{$1$};
\draw(13.63,0)node[right]{$=\mathbf{f}$};
\draw[line width=1.7pt, color=bblau, ->](0.005,1.6)--(0.005,2.5)--(13.615,2.5)--(13.615,1.6);
\draw (6.2,2.5)node[above]{\small \color{bblau}{addition}};
\draw[line width=1.7pt, color=goyel, ->](1.2,0)--(11.4,0);
\draw[line width=1.7pt, color=goyel, ->](2.3,-2.3)--(10.2,-2.3);
\draw (6,-2.3)node[above]{\small \color{goyel}{concatenation}};
\draw[line width=1.7pt, color=goyel, ->](3.6,-3.85)--(8.7,-3.85);
\draw[line width=1.7pt, color=goyel, ->](5.2,-5.025)--(6.5,-5.025);
\mirror at (0,-7.2)
\draw(0.35,-7.2)node[right]{\ldots $3\times 3$ convolutions followed by ReLU activation.};
\mirrord at (0.1,-7.45)
\draw(0.35,-7.6)node[right]{\ldots Downsampling ($2\times 2 \ \max-$pooling).};
\mirroru at (0.1,-8.1)
\draw(0.35,-8)node[right]{\ldots Upsampling followed by $3\times 3$ convolutions with ReLU as activation.};
\mirrortwo at (0,-8.4)
\draw(0.35,-8.4)node[right]{\ldots $1\times 1$ convolution followed by identity activation function.};
\draw(0.35,-8.9)node[right]{};
\draw(0,0)node[left]{$\mathbf{b}=$};
\end{tikzpicture} 
	}
    \caption{\textbf{Architecture of the residual U-net.} The number of convolution kernels (channels) is written over each layer. Long arrows indicate direct connections with subsequent concatenation or addition.}
    \label{fig:unet}
\end{figure}

The probably simplest approach is to use an explicit reconstruction function $R_\theta = \mathcal{N}_\theta \circ \A^\sharp \colon  \R^{mQ}\to\R^n$ where $\mathcal{N}_\theta$ is a neural network (NN) and $\A^\sharp$
is an operator that performs an initial reconstruction.
 In order to determine  the parameter vector  $\theta \in \R^p $ (where $p$ can be very large) that parameterizes the NN, one minimizes an error function averaged over a finite set of training data $(\mathbf{b}_k, \mathbf{f}_k)_{k=1}^N$. Here $\mathbf{f}_k$ are samples of phantoms and $\mathbf{b}_k = \A^\sharp \A(\mathbf{f}_k)$ the corresponding input images.
Then one solves the following optimization problem iteratively
\begin{equation}\label{eq:NNerror}
    \min_{\theta} \frac{1}{N}\sum_{k=1}^N \norm{\mathcal{N}_\theta(\mathbf{b}_k)-\mathbf{f}_k}_p^q \,.
\end{equation}
In particular, stochastic gradient and variants are frequently
applied to approximately minimize \eqref{eq:NNerror}.

In order to simplify the learning procedure~\cite{han2016deep},
it has been proposed to train a NN that learns the residual images $\mathbf{f} - \mathbf{b} =  \mathbf{f}  - \A^\sharp \A \mathbf{f}$. In such a situation, the reconstruction function has the form
\begin{equation} \label{eq:resnet}
    R_\theta^\text{res} = (\operatorname{Id}+\mathcal{U}_\theta)\A^\sharp \,,
\end{equation}
where $\theta \in \R^p $ is the adjustable parameter vector. A popular choice~\cite{antholzer2018deep,jin2017deep,han2016deep} for $\mathcal{U}_\theta$ is the so called U-net, which was originally designed for image segmentation~\cite{ronneberger2015unet}. The resulting NN architecture is shown in Figure~\ref{fig:unet}.  
Variants of the residual structure to increase data consistence 
have been proposed in \cite{schwab2018deep,schwab2018big}.

\subsection{Joint $\ell_1$-minimization}

In~\cite{sandbichler2018sparsification} a method based on $\ell_1$-minimization was introduced  which relies on sparsity.
An element $\mathbf{v}\in\R^n$ is called $s$-sparse, with $s \in \{1, \dots, n\}$, if it contains at most $s$ nonzero elements.
One can reconstruct $\mathbf{v}$ in a stable manner  from measurements
 $\mathbf{g} = \A\mathbf{v}$ provided that  $\A$ satisfies the restricted isometry property (RIP) of order $2s$. This property means that $(1-\delta)\norm{\z}^2 \le \norm{\A\z}^2 \le (1+\delta)\norm{z}^2$ holds for all $\z\in\R^n$ which are $2s$-sparse and the  constant
 $\delta < 1/\sqrt{2}$ \cite{CS-introduction}. Bernoulli matrices satisfy the RIP with high probability~\cite{CS-introduction}, but the subsampling matrix does not. Also it is not clear if the forward operator $\A=(\mathbf{S}\otimes\mathbf{I})\circ\mathcal{W}$ satisfies the RIP for any sampling matrix.
However the following result from inverse problems theory can still be applied in this case.\cite{grasmair2011necessary}

\noindent\textbf{Theorem 1.}
{\em 
Let $\A\in\R^{mQ\times n}$ and $\mathbf{v}\in\R^n$. Assume that the source condition holds: $
	\exists \mathbf{w}\in\R^{mQ} \colon \A^\intercal \mathbf{w}\in \sign(\mathbf{v})
	\wedge  \forall  i \in \supp(\mathbf{v}) \colon \abs{(\A^\intercal\mathbf{w})_i} < 1$, where $\sign(\mathbf{v})$ is the set-valued signum function and $\supp(\mathbf{v})$ is the set of indices with nonzero components  of $\mathbf{v}$. Further, assume that $\A$, restricted to the subspace spanned by $e_i$ for $i\in\supp(\mathbf{v})$, is injective.
    Then for any $\mathbf{g}^\delta\in\R^{mQ}$ such that $\norm{\A\mathbf{v} - \mathbf{g}^\delta}_2 \leq \delta$ and any minimizer of the $\ell_1$-Tikhonov functional,
	$\mathbf{v}_\lambda^\delta \in \argmin_{\mathbf{z}} \frac{1}{2}\norm{\A\mathbf{z}-\mathbf{g}^\delta}_2^2 + \lambda \norm{\mathbf{z}}_1$ we have  $\norm{\mathbf{v}_\lambda^\delta-\mathbf{v}}_2 = \mathcal{O}(\delta)$ provided that $\lambda \asymp \delta$.
}

It was also shown~\cite{grasmair2011necessary}, that the RIP implies the source condition of Theorem~1. Additionally, a smaller support set
$\supp(\mathbf{v})$ makes it easier to fulfill the conditions in Theorem~1. This indicates that sparsity of the unknowns is an important requirement for $\ell_1$-minimization.
The sparsity approach of \cite{sandbichler2018sparsification}
is based on the following result: \cite{sandbichler2018sparsification}
If $f$ is an initial PA source vanishing outside of $B_R$ and let $p$ be the causal solution of \eqref{eq:wave}. Then $\partial_t^2p$ is the causal solution of \eqref{eq:wave} with modified source,
$    \partial^2q(\ro, t) - \Delta_\ro q(\ro, t)
	= \delta^\prime(t) c^2\Delta\f(\ro) \quad \text{for } (\ro, t)\in \R^2\times \R_+$.
As a consequence we have that
    \begin{equation} \label{eq:cr}
	\forall \mathbf{f}\in\R^n\colon\quad \mathcal{D}_t^2\A \mathbf{f} = \A (c^2\mathcal{L}_\ro \mathbf{f})
    \end{equation}
    holds up to discretization errors, where $\mathcal{L}_\ro$ is the discrete Laplacian and $\mathcal{D}_t$ the discrete temporal  derivative.

For typical PA sources, $\mathcal{L}_\ro\mathbf{f}$ is sparse or at least compressible. Thus, based on  \eqref{eq:cr}, one  
 could first recover $\mathcal{L}_\ro$ by solving the following $\ell_1$-problem  $   \argmin_{\mathbf{z}} \{ \norm{\mathbf{z}}_1 \mid
  \A\mathbf{z}=\mathcal{D}_t^2\mathbf{g} \}$, 
and then solve the Poisson equation $\mathcal{L}_\ro \mathbf{f} = \mathbf{g}/c^2$ with zero boundary conditions in order to get $\mathbf{f}$. However, this approach leads to low frequency artifacts in the reconstructed phantom. To overcome this issue, a joint minimization approach  was introduced~\cite{sandbichler2018sparsification}, which jointly reconstructs $\mathbf{f}$ and $\mathbf{\mathcal{L}_\ro\mathbf{f}}$.
In practice, one minimizes the following
\begin{equation}\label{eq:jointminproblem}
    \min_{\mathbf{f},\mathbf{h}} \frac{1}{2}\norm{\A\mathbf{f}-\mathbf{g}}_2^2 + \frac{1}{2}\norm{\A\mathbf{h}-\mathcal{D}_t^2\mathbf{g}}_2^2 + \frac{\alpha}{2}\norm{\mathcal{L}_\ro\mathbf{f}-\mathbf{h}/c^2}_2^2 + \beta\norm{\mathbf{h}}_1 + I_C(\mathbf{f}) \,,
\end{equation}
where $\alpha$ is a tuning parameter and  $\beta$ is a regularization parameter. Moreover, $I_C$ implements a positivity constraint,
i.e. with $C=[0,\infty)^n$, the function $I_C$ is defined by 
 $I_C(\mathbf{f})=0$ if $\mathbf{f}\in C$ and $I_C(\mathbf{f})=\infty$ otherwise.

To solve \eqref{eq:jointminproblem}, one can use a proximal forward-backward splitting method~\cite{combettes2011proximal}, which is applicable to problems separable into smooth and non-smooth but convex parts. Here we  take the smooth part as $\Phi(\mathbf{f},\mathbf{h}) = 1/2\norm{\A\mathbf{f}-\mathbf{g}}_2^2 + 1/2\norm{\A\mathbf{h}-\mathcal{D}_t^2\mathbf{g}}_2^2 + \alpha/2\norm{\mathcal{L}_\ro\mathbf{f}-\mathbf{h}/c^2}_2^2$
and the non-smooth part as $\Psi(\mathbf{f},\mathbf{h}) =  \beta\norm{\mathbf{h}}_1 + I_C(\mathbf{f})$.
We need to calculate the proximity operator of the non-smooth parts,
which, for a convex function $F\colon\R^n\to\R$, is defined  by
$    \prox_F(\mathbf{f}) \triangleq \argmin \{ F(\mathbf{z}) + \tfrac{1}{2}\lVert\mathbf{f}
    -\mathbf{z}\rVert_2^2 \mid \mathbf{z}\in\R^n \}$. In our case, the proximity operator can be computed explicitly and component-wise,
    $\prox_{\Psi} (\mathbf{f},\mathbf{h}) = [\prox_{I_C}(\mathbf{f}), \prox_{\beta\norm{\edot}_1}(\mathbf{h})]$. 
On the other hand, the gradients of $\Phi$ can also be calculated 
explicitly. The resulting joint $\ell^1$-minimization algorithm for   \eqref{eq:jointminproblem} is given by the  following iterative scheme
\begin{equation}\label{eq:jointalgorithm}
    \begin{aligned}
	\mathbf{f}^{k+1} =& \prox_{I_C}\left( \mathbf{f}^k - \mu 
	\left( \A^\intercal(\A\mathbf{f}^k-\mathbf{g}) - \alpha \mathcal{L}(\mathcal{L}\mathbf{f}^k - \mathbf{h}^k/c^2) \right) 
	\right)\\
	\mathbf{h}^{k+1} =& \prox_{\beta\norm{\edot}_1}\left( \mathbf{h}^k - \mu \left( \A^\intercal(\A\mathbf{h}^k-\mathcal{D}_t^2\mathbf{g}) - \frac{\alpha}{c^2}(\mathcal{L}_\ro\mathbf{f}^k-\mathbf{h}^k/c^2) \right)  \right)\,,
    \end{aligned}
\end{equation}
with starting points $\mathbf{f}^0=\mathbf{h}^0=\mathbf{0}$ and 
$ \prox_{I_C}(\mathbf{f}) = (\max(\mathbf{f}_i,0))_i$,
$ \prox_{\beta\norm{\edot}_1}(\mathbf{h}) = (\max(\abs{\mathbf{h}_i}-\beta,0)\sign(\mathbf{h}_i)) $.

\section{Nett: variational regularization with neural networks}
\label{sec:nett}

Standard  deep learning  approaches have the disadvantage that
they may perform badly on images that are very different from the ones included in the training set. In order to address this issue, iterative networks~\cite{adler2018learned,kobler2017variational,kelly2017deep} or
data invariant regularizing networks \cite{schwab2018big} have been proposed. Another method to enhance data consistency, which we will use in this paper, is based on generalized Tikhonov regularization using a
learned regularization term~\cite{li2018nett}.

\subsection{NETT framework}

The basic idea is to consider minimizers  of the
unconstrained  optimization problem
\begin{equation}\label{eq:nett}
    \min_{\mathbf{f}} \frac{1}{2}\norm{\A\mathbf{f}-\mathbf{g}}_2^2 + \frac{\lambda}{2} \, \mathcal{R}(\mathbf{f}) \,,
\end{equation}
where   $\mathcal{R} \colon  \R^n \to [0, \infty]$ is a trained
regularizer and $\lambda > 0$ is the regularization parameter.
The  resulting reconstruction approach  is called NETT (for network Tikhonov regularization), as it is a generalized  form of Tikhonov regularization using a NN as trained regularizer.

In \cite{li2018nett} it  has been shown that under reasonable conditions,
the  NETT approach is well-posed and yields a convergent regularization method. In particular, minimizers   of   \eqref{eq:nett} exist,
are stable with respect to data perturbations, and minimizers of \eqref{eq:nett} converge to   $\mathcal{R}$-minimizing solutions of
the equation $\A\mathbf{f} = \mathbf{g}$ as the noise  level
goes to zero.
For the regularizer proposed in  \cite{li2018nett}  these conditions  are difficult to be verified.
  We propose a variation of the trained regularizer 
  using a simpler network architecture and  different training strategy.

\subsection{Construction of regularizers}
\label{sec:nett-train}

\begin{figure}[thb!]
    \centering
    \resizebox{!}{0.2\textwidth}{
	\begin{tikzpicture}
\filldraw[color=black](0,1.5)--(0,-1.5)--(0.01,-1.5)--(0.01,1.5)--cycle;
\draw(0.005,1.5)node[above]{$1$};
\mirror at (0.1,0)
\filldraw[color=black](0.45,1.5)--(0.45,-1.5)--(0.55,-1.5)--(0.55,1.5)--cycle;
\draw(0.5,1.5)node[above]{$64$};
\mirror at (0.6,0)
\filldraw[color=black](0.95,1.5)--(0.95,-1.5)--(1.1,-1.5)--(1.1,1.5)--cycle;
\draw(1.075,1.5)node[above]{$32$};
\mirrortwo at (1.2,0)
\filldraw[color=black](1.55,1.5)--(1.55,-1.5)--(1.56,-1.5)--(1.56,1.5)--cycle;
\draw(1.55,1.5)node[above]{$1$};

\draw(0,0)node[left]{$\mathbf{f}=$};

\draw(1.57,0)node[right]{$=\mathbf{x} \to \lVert \cdot \lVert^2_2$};
\end{tikzpicture}
	}
    \caption{\textbf{Network structure for the trained  regularizer.}
   The first  two convolutional layers use  $3\times 3$ convolutions followed by ReLU activations. The last convolutional layer  (green arrow) uses
    $3\times 3$ convolutions  not followed by an activation function.}
    \label{fig:regnet}
\end{figure}
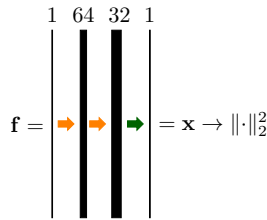

For the  regularizer in \eqref{eq:nett} we make the ansatz
\begin{equation} \label{eq:nreg}
	\mathcal{R}(\mathbf{f})
	= \norm{\mathcal{V}_\theta (\mathbf{f})}_F^2 \,,
\end{equation}
where $\norm{\cdot}_F$ is the Frobenius norm and
$\mathcal{V}_\theta\colon\R^n \to \R^n$ is a trained
NN.   In this work, we  use the simple  NN architecture  shown in
Figure~\ref{fig:regnet}, that consist of three convolutional layers.
Clearly, one could also use more complicated   network structures.

To  train the network $\mathcal{V}_\theta$  we choose a set of phantoms
$(\mathbf{f}_k)_{k=1}^{N_1 + N_2}$  and compute initial reconstructions
$\mathbf{b}_k  =  \A^\sharp \A\mathbf{f}_k$.
We then  define a  training set of input/output pairs
$(\mathbf{x}_k, \mathbf{y}_k)_{k=1}^{N_1+N_2}$
in the following way:
\begin{equation}\label{eq:regdata}
\begin{aligned}
    \mathbf{x}_k &= \mathbf{b}_k\,,
    && \mathbf{y}_k = \mathbf{b}_k - \mathbf{f}_k \quad &\text{for } k&=1,\ldots,N_1\\
    \mathbf{x}_{k} &= \mathbf{f}_{k} \,,
    && \mathbf{y}_{k} = \mathbf{0} &\text{for } k& = N_1+1,\ldots,N_1+N_2 \,.
\end{aligned}
\end{equation}
The parameters in the network $ \mathcal{V}_\theta$ are optimized 
to approximately map $\mathbf{x}_k$ to $ \mathbf{y}_k$.
For that purpose we minimize the mean absolute error
\begin{equation}\label{eq:nettmin}
    E(\theta)
    = \frac{1}{N_1 + N_2} \sum_{k=1}^{N_1+N_2}\norm{\mathcal{V}_\theta (\mathbf{x}_k)-\mathbf{y}_k}_1
\,.
\end{equation}
averaged of the training set.

Note that the trained regularizer depends on the forward operator 
as well as on the initial reconstruction operator.
If the equation $\A \x = \y$ is undetermined, it is reasonable to
take the initial reconstruction operator as  
right inverse $\A^\sharp \colon \R^{mQ} \to \R^n$, i.e. 
$\A \A^\sharp \y = \y $ holds for exact data $\y$. The residual images
$\A \A^\sharp \x - \x$  in this case are contained in the
null-space of the forward operator and  trained regularizer 
finds and penalizes the component of $\x$  in the null space.
This reveals connections of NETT with the null space 
approach of
\cite{schwab2018deep}. However, for training the regularizer one can  consider other  initial reconstructions  that add undesirable structures. 
For CS-PAT, the forward operator has the form 
$\A = \mathcal{S}\circ\mathcal{W}$  and we use  
$\A^\sharp$ as a numerical approximation of
$\mathbf{W}^\intercal \circ \mathbf{S}^\intercal$.

For comparison purpose we also test \eqref{eq:nett} with the
deterministic regularizer
$\mathcal{R}(\mathbf{f}) = \norm{D_x\mathbf{f}}_F^2 + \norm{D_y\mathbf{f}}_F^2$, where $D_x$ and $D_y$ 
denote the discrete derivatives in $x$ and $y$ direction, respectively. We can write these derivatives
as convolution operations with the convolution  kernels
$ k_x = 2^{-1/2}[-1, 1]$  and $k_y = 2^{-1/2}[-1,1]^\intercal$.
The deterministic regularizer therefore has the form
$\norm{\mathcal{D}(\mathbf{f})}^2_F $, where $\mathcal{D}=[D_x; D_y]$
is a convolutional NN with no hidden layer, a two-channel output layer, and
fixed non-adjustable parameters. We therefore call  the
 resulting  reconstruction  approach the deterministic NETT.
 Notice that the  deterministic NETT  is equal to standard
 $H^1$-regularization.

\subsection{Minimization of the NETT functional}

Once the weights of the network have been trained,
we reconstruct the phantom by  iteratively minimizing
the NETT functional \eqref{eq:nett}.
For that  purpose we use the following incremental
gradient  algorithm:
\begin{equation}\label{eq:splitgradient}
    \begin{aligned}
	\hat{\mathbf{f}}^{k+1} &= \mathbf{f}^k - \mu \left( \A^T(\A\mathbf{f}^k - \mathbf{g})\right)\\
	\mathbf{f}^{k+1} &= \hat{\mathbf{f}}^{k} - \mu\lambda \left( \nabla_\mathbf{f} \mathcal{V}_\theta(\mathbf{f}^k)\right)\,.
    \end{aligned}
\end{equation}
Note that the derivative $\nabla_\mathbf{f}\mathcal{V}_\theta$ is with respect to the input of the NN and not its parameters. This gradient
can be calculated by standard deep learning software.
Note that by fixing the number of iterations, the iteration
\eqref{eq:splitgradient} shares some similarities with iterative and variational networks~\cite{adler2018learned,kobler2017variational}.

The NETT convergence theory of \cite{li2018nett} requires
the regularization functional $\mathcal{R}$  to be proper,
coercive and  weakly lower semi-continuous.
In the finite-dimensional setting considered above, this is equivalent
to the  coercivity  condition
$\forall  \mathbf{f} \in \R^n \colon \norm{\mathcal{V}_\theta (\mathbf{f})}_F^2  \geq c \norm{\mathbf{f}}_F^2$. 
We did not explicitly account for this
condition in the network construction. To enforce stability, we
may combine  \eqref{eq:splitgradient} with early stopping.
As an alternative strategy, we might adjust the training process,
or replace the trained regularizer by  $ \norm{\mathcal{V}_\theta (\mathbf{f})}_F^2 + a \norm{ \mathbf{f}}_F^2$ with some constant $a>0$.

\section{Numerical results}
\label{sec:numerical}

In this section we  present  numerical results including
simulated as well as experimental data in order to compare the methods introduced  in the previous section.

\subsection{Implementation details}

We use Keras~\cite{keras} and Tensorflow~\cite{tensorflow} for implementation and optimization of the NNs. The FBP algorithm and the joint $\ell_1$-algorithm are  implemented in MATLAB.
We ran all our experiments on a computer using an Intel i7-6850K and an NVIDIA 1080Ti.

Any discrete PA source $\mathbf{f}\in\R^n$ with $n=128^2$ consists of discrete samples of the continuous source at a $128\times 128$ grid covering the square $[\SI{-5}{\micro\metre}, \SI{9}{\micro\metre}]\times[\SI{-12.5}{\micro\metre}, \SI{1.5}{\micro\metre}]$. The full wave data $\mathbf{g}\in\R^{MQ}$ correspond to $M=240$ sensor locations on the circle of radius $\SI{40}{\micro\metre}$ and polar angles in the interval $[35^\circ, 324^\circ]$ and $Q=747$ equidistant temporal samples in $[0,T]$ with $cT=\SI{4.9749d-2}{\micro\metre}$. The sound speed is taken as $c = \SI{1.4907d3}{\metre\per\second}$.

The wave data are  simulated by discretization of the wave equation, and~\eqref{eq:fbp} is implemented using the standard FBP approach~\cite{FinHalRak07}. This gives us a forward operator $\mathcal{W}\colon\R^n\to\R^{mQ}$ and an unmatched adjoint $\mathcal{B}\colon\R^{nQ} \to \R^n$.
We consider $m=60$ CS measurements, which yields a compression ratio of $4$. We also generated a noisy dataset by adding 7\% Gaussian white noise to the measurement data.
For the sampling matrices $\mathbf{S}\in\R^{m\times M}$ we use the deterministic sparse subsampling matrix with entries
$\mathbf{S}[i,j] =2$ if $j=4(i-1)+1$ and $0$ otherwise,
and the random Bernoulli matrix where each entry is taken independently at random as $\pm 1/\sqrt{m}$ with equal probability.

\subsection{Reconstruction methods}

For the two-step as well as the NETT approach  we use
$\A^\sharp  \coloneqq  \mathcal{B} \circ \mathbf{S}^\intercal$ as initial reconstruction. To train the networks we generate a dataset of $N=500$ training examples $(\mathbf{b}_k,\mathbf{f}_k)_{k=1}^N$ where $\mathbf{f}_k$ are taken as projection images from three-dimensional blood vessel data as described in~\cite{schwab2018dalnet}, and  
$\mathbf{b}_k = \A^\sharp \A\mathbf{f}_k$.
To train the residual U-net we minimize the  mean absolute
error $\frac{1}{N} \sum_{k=1}^N\norm{(\operatorname{Id}+\mathcal{U}_\theta)(\mathbf{b}_k) - \mathbf{f}_k}_1$.
For  NETT regularization, we use training   data as  in \eqref{eq:regdata}
with $N_1=N_2 = N$  and $\mathbf{f}_{k+N}
= \mathbf{f}_k$ and minimize \eqref{eq:nettmin}.
In both cases we use the Adam optimizer~\cite{kingma2014adam} with 300 epochs and a learning rate of 0.0005.

For the joint $\ell_1$-method we use 70 iterations of \eqref{eq:jointalgorithm} with the parameters $\alpha=0.001,\, \beta=0.005$ and $\mu=0.125$ ($\mu=0.03125$ for noisy data) for Bernoulli sampling, and $\mu=0.0625$ ($\mu=0.03125$ for noisy data) for sparse sampling.
For  NETT regularization we use 10 iterations of \eqref{eq:splitgradient} with $\mu=0.5$ (0.7 for noisy phantoms) and $\lambda=0.5$ for the deterministic regularization network we use $\mu=0.5$ and $\lambda=0.35$.
All hyper-parameters have been selected by hand to get good visual results and no hyper-parameter optimization has been performed.

\begin{figure}[thb!]
    \includegraphics[width=\textwidth]{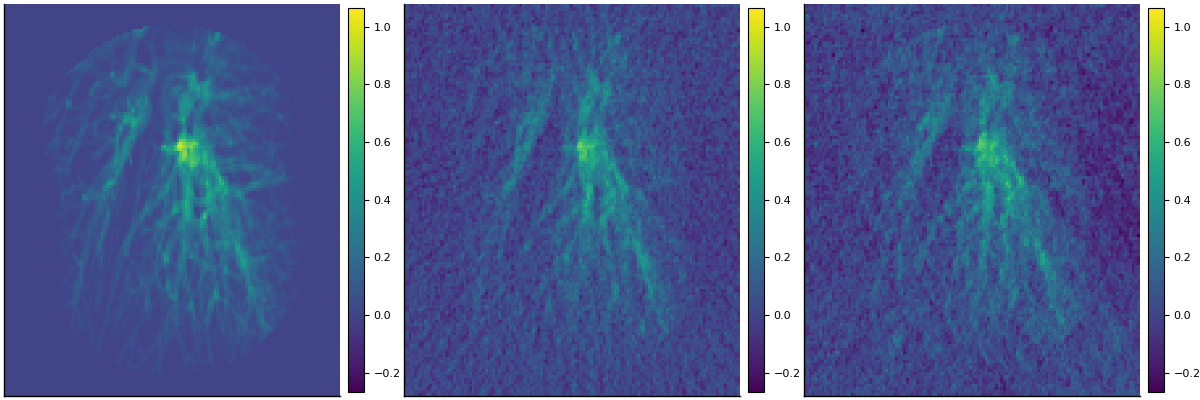}
    \caption{\textbf{Sample from the blood-vessel data set.}
    Left: Ground truth phantom.
    Middle: Initial reconstruction using $\A^\sharp$ from sparse data.
    Right: Initial reconstruction using $\A^\sharp$ from Bernoulli data.}
    \label{fig:blood1}
\end{figure}

\begin{figure}[thb!]
    \includegraphics[width=\textwidth]{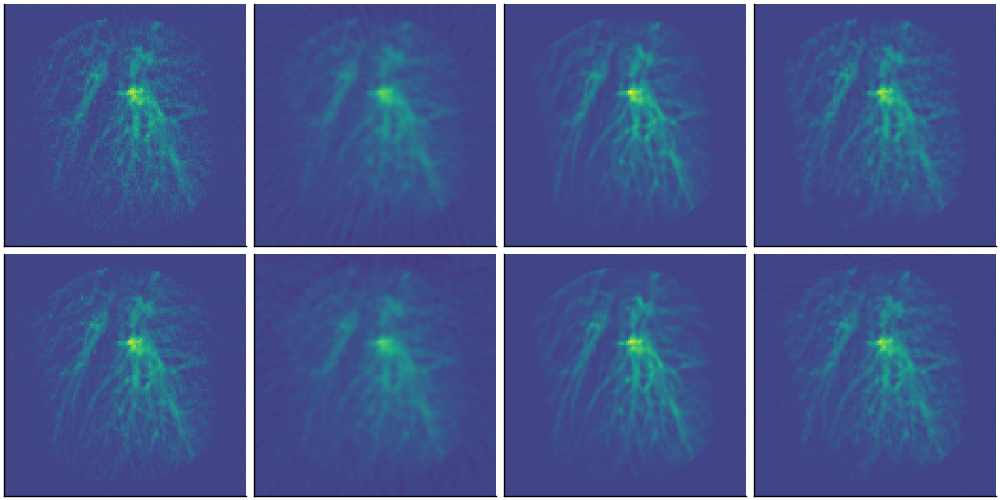}
    \caption{\textbf{Reconstruction results for using simulated data.}
	Top Row: Reconstructions from sparse data.
    Bottom Row:    Reconstructions from Bernoulli data.
    First Column:  Joint $\ell_1$-algorithm.
    Second Column: $H^1$-regularization (deterministic NETT).
    Third Column:  Residual U-Net.
    Fourth Column: NETT.
    }
    \label{fig:blood2}
\end{figure}

\begin{figure}[thb!]
    \includegraphics[width=\textwidth]{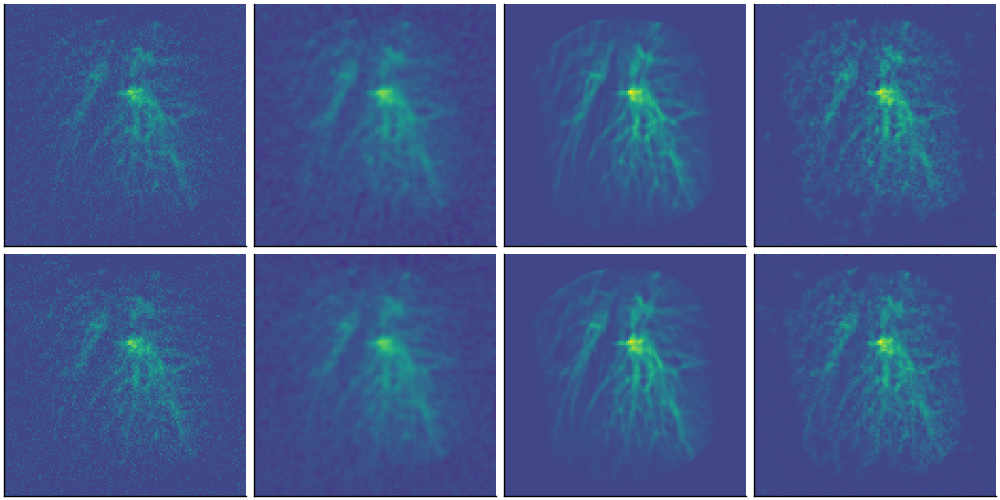}
    \caption{\textbf{Reconstruction results  using noisy data.}
	Top Row: Reconstructions from sparse data.
    Bottom Row:    Reconstructions from Bernoulli data.
    First Column:  Joint $\ell_1$-algorithm.
    Second Column: $H^1$-regularization (deterministic NETT).
    Third Column:  Residual U-Net.
    Fourth Column: NETT.
    }
    \label{fig:blood2-noise}
\end{figure}

\begin{table}[thb!]\centering
    \begin{subtable}{0.47\textwidth}
    \begin{tabular}{|r|c|c|c|c|}
\toprule
	Method & MSE & RMAE & PSNR & SSIM\\
	\midrule
	FBP & 0.00371 & 4.52 & 24.92 & 0.39\\
	$\ell_1$ & 0.00075 & 1.71 & 31.69 & 0.76\\
	$H^1$ & 0.00108 & 2.02 & 30.05 & 0.75\\
	U-net & 0.00072 & 1.46 & 31.94 & 0.86\\
	NETT & \hlc{0.00048} & \hlc{1.39} & \hlc{33.56} & \hlc{0.89}\\
	\bottomrule
    \end{tabular}
	\vspace{0.5em}
    \end{subtable}
    \hspace{1em}
    \begin{subtable}{0.47\textwidth}
	\begin{tabular}{|r|c|c|c|c|}
	\toprule
	Method & MSE & RMAE & PSNR & SSIM\\
	\midrule
	FBP & 0.00472 & 5.38 & 23.44 & 0.32\\
	$\ell_1$ & \hlc{0.00028} & \hlc{1.08} & \hlc{35.79} & 0.86\\
	$H^1$ & 0.00126 & 2.22 & 29.27 & 0.7\\
	U-net & 0.00096 & 1.59 & 31.37 & 0.85\\
	NETT & 0.00045 & 1.43 & 33.68 & \hlc{0.88}\\
    \bottomrule
    \end{tabular}
	\vspace{0.5em}
    \end{subtable}
    \caption{\textbf{Averaged performance for noise-free data.}
    Left: Sparse sampling. Right: Bernoulli sampling. Best values are highlighted.}
    \label{tab:meanerror}
\end{table}

\begin{table}[thb!]\centering
    \begin{subtable}{0.47\textwidth}
    \begin{tabular}{|r|c|c|c|c|}
    \toprule
	Method & MSE & RMAE & PSNR & SSIM\\
	\midrule
	FBP & 0.00371 & 4.52 & 24.92 & 0.39\\
	$\ell_1$ &0.00245 & 3.24 & 26.23 & 0.46\\
	$H^1$ & 0.00134 & 2.54 & 28.95 & 0.62\\
	U-net & \hlc{0.00076} & \hlc{1.52} & \hlc{31.61} & \hlc{0.85}\\
	NETT & 0.0012 & 2.24 & 29.47 & 0.75\\
    \bottomrule
    \end{tabular}
	\vspace{0.5em}
    \end{subtable}
        \hspace{1em}
    \begin{subtable}{0.47\textwidth}
	\begin{tabular}{|r|c|c|c|c|}
	\toprule
	Method & MSE & RMAE & PSNR & SSIM\\
	\midrule
	FBP & 0.00474 & 5.38 & 23.44 & 0.32\\
	$\ell_1$ & 0.00277 & 3.47 & 25.68 & 0.43\\
	$H^1$ & 0.00139 & 2.47 & 28.77 & 0.65\\
	U-net & 0.00098 & \hlc{1.61} & \hlc{31.19} & \hlc{0.84}\\
	NETT & \hlc{0.0008} & 1.82 & 31.17 & 0.82\\
    \bottomrule
    \end{tabular}
	\vspace{0.5em}
    \end{subtable}
    \caption{\textbf{Averaged performance for data including 7\% noise.}
    Left: Sparse sampling. Right: Bernoulli sampling. Best values are highlighted.}
    \label{tab:meanerror-noise}
\end{table}

\subsection{Results for simulated data}

For the offline evaluation, we investigate  performance on 10 blood vessel phantoms (not contained in the training set). We consider sparse sampling and Bernoulli sampling. One of the evaluation phantoms
 and reconstruction results for noise-free and for noisy data are shown  in Figures~\ref{fig:blood1}-\ref{fig:blood2-noise}.

To quantitatively evaluate the reconstruction quality, we calculated the relative mean absolute error (RMAE), the mean squared error (MSE), the peak signal to noise ratio (PSNR) and the structured similarity index (SSIM). The performance measures 
averaged over the 10 evaluation phantoms
are shown in Table~\ref{tab:meanerror} for noise-free data,
and in Table~\ref{tab:meanerror-noise} for the noisy data case.
We can see that the learned approaches work particular well for  the sparse sampling matrix.
The residual U-net seems to be better than NETT at denoising the image, which results in lower errors for noisy data. We further observe that our simple trained regularizer in any case
performs better than the (very simple) deterministic one.
The $\ell_1$-minimization approach works best for noise-free Bernoulli measurements, but it is outperformed by the learned approaches in the noisy data case.
Since $\ell_1$-minimization needs many iterations we did not iterate until convergence;
for the presented reconstructions it is  already one order of magnitude slower than the other methods.
In the noisy data case, we observe that the NETT yield smaller MSE than the residual U-net  for Bernoulli measurements, while the U-Net approach works best for the sparse measurements.

\subsection{Results for experimental data}

We used experimental data measured by the PAT device using integrating line sensors as described in~\cite{bauer2017all,sandbichler2018sparsification}.
The setup corresponds to the one used for our simulated phantoms.
The first sample is a simple cross phantom and the second is a leaf phantom with many fine structures. 
We only test sparse measurements, since the current experimental
setup does not support Bernoulli measurements. For the residual U-net and the  NETT with  trained  regularizer we use the networks trained on the blood vessel data set
as described above.

\begin{figure}[thb!]\centering
    \includegraphics[width=\textwidth]{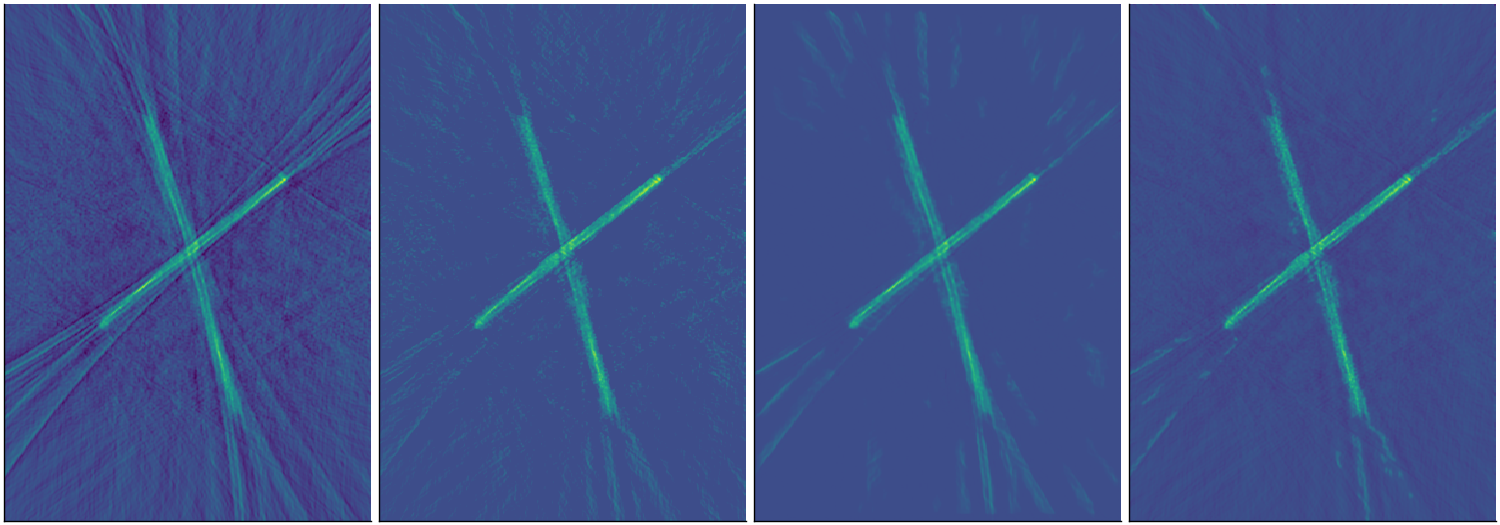}
    \caption{\textbf{Reconstructions of the cross phantom from experimental data.}
     The reconstruction results are obtained for the sparse sampling pattern
     using the following algorithms, from left to right:
     FBP, $\ell_1$-minimization, U-net and NETT.}
    \label{fig:cross}
\end{figure}

\begin{figure}[thb!]\centering
    \includegraphics[width= \textwidth]{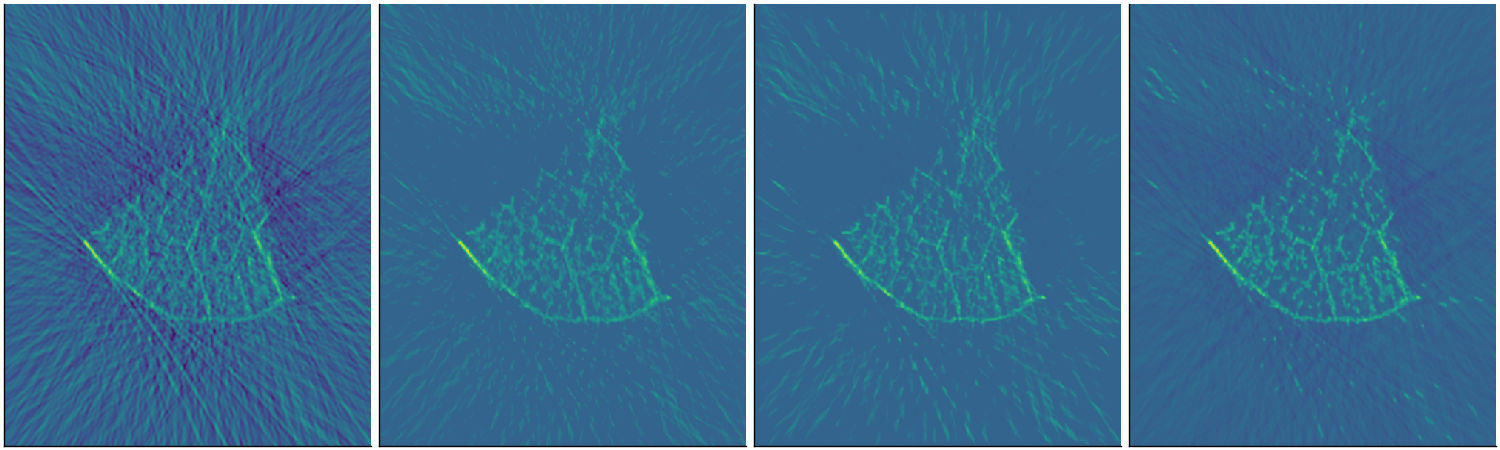}
    \caption{\textbf{Reconstructions of the leaf phantom from experimental data.}
     The reconstruction results are obtained for the sparse sampling pattern
     using the following algorithms, from left to right:
     FBP, $\ell_1$-minimization, U-net and NETT.}
    \label{fig:leaf}
\end{figure} 
Reconstruction results  for the cross phantom are shown in Figure~\ref{fig:cross}.
All methods yield quite accurate results and all structures are resolved  clearly.
However, none of the methods is able to completely remove the strong artifacts extending
from the corners of the cross phantom. The U-net seems to do the best job
by removing nearly all artifacts.
Reconstruction results  for the leaf phantom  are shown in Figure~\ref{fig:leaf}.
Reconstructing the leaf phantom is more challenging, since it contains very fine structures.
None of the methods is able to resolve them completely.
The artifacts outside of the object are present in all
reconstructed images. However, the reconstruction using joint $\ell^1$-minimization and the deep-learning approaches (NETT and residual U-net) yield satisfying results. Nevertheless, future work is required to improve the reconstruction quality for real-word data.

\section{Discussion and conclusion}

We studied deep-learning approaches (NETT and the residual U-net)
and joint $\ell_1$-minimization  for CS-PAT using either sparse sampling 
or Bernoulli measurements.
All methods work  well for both types of measurements. For exact data, iterative approaches,  which use the forward operator in each step, work better than the residual U-net for Bernoulli measurements. Incorporating different sampling strategies directly in the experimental setup is an interesting line for future research. NN based algorithms perform well with noisy data, but still are open for  improvements for experimental data.
This suggests that our simulated training data is different from the measured real-world data. Developing more accurate forward models and improving training data are an important future goal.
The learned regularizer has a quite  simple network structure. Investigating more complex network architectures and modified training strategies will be investigated in future work.

\section*{Acknowledgments}
S.A., J.B and M.H.  acknowledge support of the Austrian Science Fund (FWF), project P 30747-N32.

\end{document}